\documentclass[twoside,12pt]{article}

\usepackage{amsmath}
\usepackage{amssymb}
\usepackage{amscd}
\usepackage{amsthm}

\numberwithin{equation}{section}

\theoremstyle{plain}

\theoremstyle{remark}

\theoremstyle{definition}

\newcommand{\Z}{\mathbb Z}
\newcommand{\R}{\mathbb R}
\newcommand{\C}{\mathbb C}

\def\ra{\rightarrow}

\def\e{\emph}
\def\i{\infty}

\def\b{\begin}
\newcommand{\ol}{\overline}

\begin{document}

\title{
{Rigidity  of  quasiconformal maps on Carnot groups}}
\author{Xiangdong Xie\footnote{Partially supported by NSF grant DMS--1265735.}}
\date{  }

\maketitle

\begin{abstract}
We show that quasiconformal maps on many Carnot groups must be biLipschitz. 
 In particular, this is the case for  
 $2$-step Carnot groups  with  reducible first layer.  These results have implications for the   rigidity of quasiisometries between negatively curved solvable Lie groups.

\end{abstract}




 {\small {\bf{Mathematics Subject
Classification (2000).}}  22E25,   30L10,   53C17. 






\setcounter{section}{0} 
  \setcounter{subsection}{0}

\section{Introduction}\label{s0}

In this paper we  establish rigidity results for quasiconformal maps between Carnot groups.   We show that 
quasiconformal maps on many Carnot groups must be biLipschitz.

Let $\mathcal N=V_1\oplus \cdots \oplus V_r$ be a Carnot   Lie algebra. 
  For each $t>0$, the dilation  $\lambda_t: \mathcal N\ra \mathcal N$ is defined 
 by $\lambda_t(v)=t^i v$ for $v\in V_i$.
 An isomorphism $A:  \mathcal N\ra \mathcal N$ is  called a {\it{graded isomorphism}}
   if $A$ commutes with the dilations $\lambda_t$ for all $t>0$; that is, 
 $A\circ \lambda_t=\lambda_t\circ A$ for all $t>0$.    Let $\text{Aut}_g(\mathcal N)$ be the group of graded isomorphisms of $\mathcal N$. We say $V_1$ is 
{\it reducible} (or the first layer of $\mathcal N$ is reducible) if  there is a non-trivial proper
 linear subspace $W_1\subset V_1$ such that $A(W_1)=W_1$ for every 
 $A\in \text{Aut}_g(\mathcal N)$. 


Let $K\ge 1$ and $C>0$. A bijection $F:X\ra Y$ between two
     metric spaces is called a $(K, C)$-\e{quasi-similarity}  if
\[
   \frac{C}{K}\, d(x,y)\le d(F(x), F(y))\le C\,K\, d(x,y)
\]
for all $x,y \in X$.

Clearly a map is a quasi-similarity if and only if it is biLipschitz.  
 The point here is that often there is control on $K$  but not on $C$.  In this case, the notion of quasi-similarity provides more information about the distortion.

We say that a Carnot group  $N$   is {\it quasisymmetrically rigid} if  
 every $\eta$-quasisymmetric map  is a $(K, C)$-quasi-similarity,  where $K$ is a constant depending only on $\eta$.   See Section \ref{maps} for the definition of  quasisymmetric map.

Here is our first result:

\b{Th}\label{main1}
  Let $N$ be a $2$-step Carnot group. If the first layer of its Lie algebra is reducible,
  then $N$ is quasisymmetrically rigid. 
\end{Th}

Of course, for many Carnot Lie algebras, the first layer is not reducible. On the other hand,  every Carnot algebra is a subalgebra of some Carnot algebra with reducible first layer.  See   Section \ref{product}   for some constructions. 




For general Carnot groups, we have the following:  

\b{Th}\label{main2}
Suppose $W_1\subset V_1$ is a non-trivial proper subspace that is invariant under the action of 
 $\text{Aut}_g(\mathcal{N})$. 
  If   there is some $X\in V_1\backslash W_1$ such that $[X, W_1]\subset [W_1, W_1]$,  
  then $N$ is quasisymmetrically rigid.
\end{Th}

 The following result is very useful for induction argument:

\b{Th}\label{main3}
 Suppose $W_1\subset V_1$ is a non-trivial proper subspace that is invariant under the action of 
 $\text{Aut}_g(\mathcal{N})$.   Let $W$ be   the connected subgroup of $N$ with  Lie algebra $<W_1>$.      If $W$ is quasisymmetrically rigid, then so is 
$N$.  
\end{Th}

   The only Carnot groups known to admit non-biLipschitz quasiconformal maps are the Heisenberg groups  \cite{B}.    We conjecture that the Heisenberg groups are the only exceptions.

The results in this paper extend previous results by the author 
 for reducible Carnot groups  \cite{X3}  and model Filiform groups   \cite{X2}. 
  In a forthcoming paper \cite{X1}  the author proves that even 
most non-rigid Carnot groups are quasisymmetrically rigid.

The first rigidity theorem about quasiconformal maps of Carnot groups is due to  Pansu.  He proved  that (\cite{P})   every quasiconformal map 
  of the quarternionic Heisenberg group is a composition of left translations and graded automorphisms. 
   In particular, the space of quasiconformal maps
 is finite dimensional.   Our results are of somehow different nature.   Our results imply that many Carnot groups with rank  at  most  one 
 are quasisymmetrically rigid. It will be shown in \cite{X1} that the space of biLipschitz maps of Carnot groups with rank   at   most one 
 is infinite dimensional.  See Section \ref{rank} for definition of rank and Section \ref{examples} for examples of  
quasisymmetrically rigid Carnot groups with rank at most one. 


  The results in this paper have implications for the large scale geometry of
 negatively curved 
 homogeneous manifolds. 
  Each Carnot group arises  as  the (one point complement  of) ideal boundary  of some negatively
  curved   homogeneous manifold  \cite{H}.  
   Furthermore,  each quasiisometry of the negatively curved homogeneous manifold associated to a  quasisymmetrically  rigid Carnot group  
   is a rough isometry, that is, it must preserve the distance up to an additive constant.     

In Section \ref{prelimi}  we recall basic definitions and facts about Carnot groups. 
  In  Section \ref{QCBL}  we prove Theorems \ref{main2} and \ref{main3}.  
 In Section \ref{invariant}   we  explain how to find invariant subspaces and 
  introduce   constructions that produce Carnot algebras with  reducible first layer
 (so that the Theorems in this paper  can be applied).  In Section \ref{2stepcase}  we  prove Theorem \ref{main1}. Finally in Section \ref{examples}   we 
   apply our results to  3 explicit examples.

\noindent {\bf{Acknowledgment}}. {This work was completed while the author was attending the workshop 
\lq\lq Interactions between analysis and geometry" at IPAM,  University of California at Los Angeles  from March to June 2013.  I would like to thank IPAM  for   financial support, excellent working conditions and  
   conducive atmosphere.
 I also would like to thank
   Tullia Dymarz   and 
 David Freeman  for discussions about Carnot groups.  }

\section{Preliminaries}\label{prelimi}

In this Section we collect definitions and results that shall be needed later. 
  We first recall the basic definitions related to Carnot groups  in Subsection \ref{basics}.
   Then we  review      the BCH formula (Subsection \ref{BCH}),
 the definition  of   quasisymmetric maps  (Subsection \ref{maps})  and    
  Pansu differentiability theorem (Subsection \ref{pansud}). 

\subsection{Carnot algebras and Carnot groups}\label{basics}

A   \e{Carnot Lie algebra} is a finite dimensional Lie algebra  
$\mathcal G$  over $\R$      together with  a direct sum   decomposition  
    $\mathcal
G=V_1\oplus V_2\oplus\cdots \oplus V_r$ 
  of   non-trivial   vector subspaces
 such that $[V_1,
V_i]=V_{i+1}$ for all $1\le i\le r$,
    where we set $V_{r+1}=\{0\}$.  The integer $r$ is called the
    degree of nilpotency of $\mathcal
G$. Every Carnot algebra
 $\mathcal
G=V_1\oplus V_2\oplus\cdots \oplus V_r$  admits a one-parameter
 family of automorphisms $\lambda_t: \mathcal
G \ra \mathcal G$, $t\in (0, \i)$, where
 $\lambda_t(x)=t^i x$ for  $x\in V_i$.
  Let   $\mathcal
G=V_1\oplus V_2\oplus\cdots \oplus V_r$
    and $\mathcal
G'=V'_1\oplus V'_2\oplus\cdots \oplus V'_s$  be two  Carnot
    algebras.
  A Lie algebra  homomorphism
  $\phi: \mathcal
G\ra \mathcal G'$
     is graded if $\phi$ commutes with $\lambda_t$ for
  all $t>0$; that is, if $\phi\circ \lambda_t=\lambda_t\circ
  \phi$.  We observe that $\phi(V_i)\subset V'_i$ for all $1\le i\le
  r$.

A  simply connected nilpotent Lie group is a  \e{Carnot group}
 if its Lie algebra is a Carnot algebra.
    Let $G$ be a Carnot group with Lie algebra
      $\mathcal G=V_1\oplus \cdots \oplus V_r$.  The subspace $V_1$ defines
      a left invariant distribution  $H G\subset TG$ on $G$.    We fix a left invariant inner product on
          $HG$.
           An
      absolutely continuous curve $\gamma$ in $G$  whose velocity vector
       $\gamma'(t)$  is contained in  $H_{\gamma(t)} G$ for a.e. $t$
        is called  a horizontal curve.
          By Chow's theorem ([BR, Theorem 2.4]),   any two points
  of $G$ can be  connected by horizontal curves. Let $p,q\in G$, the
  \e{Carnot   metric} $d_c(p,q)$   between them is defined as
  the infimum of length of horizontal curves that join $p$ and $q$.

  Since the inner product on   $HG$ is left invariant, the Carnot
  metric on $G$ is also left invariant.  Different choices of inner
  product on $HG$ result in Carnot metrics that are biLipschitz
  equivalent.
    The Hausdorff dimension of $G$ with respect to  a  Carnot metric
    is  given by $\sum_{i=1}^r i\cdot \dim(V_i)$.

Recall that, for a simply connected nilpotent Lie group $G$ with Lie
algebra $\mathcal G$, the exponential map
  $\exp: {\mathcal G}\ra G$ is a diffeomorphism.  Under this identification the Lesbegue   measure on 
  $\mathcal G$  is a  Haar measure on $G$.  
 Furthermore, the
  exponential map induces 
  a  one-to-one correspondence between
    Lie subalgebras of $\mathcal G$   and
  connected Lie subgroups of $G$.

Let $G$ be a Carnot group with Lie algebra
      $\mathcal G=V_1\oplus \cdots \oplus  V_r$.
         Since $\lambda_t: \mathcal G\ra  \mathcal G$ ($t>0$) is a Lie algebra
         automorphism  and $G$ is simply connected,  there is a
         unique  Lie group automorphism  $\Lambda_t: G\ra G$ whose
         differential  at the identity is $\lambda_t$.
       For each $t>0$,
      $\Lambda_t$ is a
             similarity with respect to the Carnot metric:  $d(\Lambda_t(p),
             \Lambda_t(q))=t\, d(p,q)$ for any two points $p, q\in
             G$.  A  Lie group homomorphism
                $f: G\ra G'$ between two Carnot groups is a
                  graded homomorphism if  it commutes with
                    $\Lambda_t$ for all  $t>0$; that is, if
                      $f\circ \Lambda_t=\Lambda_t\circ  f$.
                         Notice that,
   a Lie group homomorphism
                $f: G\ra G'$ between two Carnot groups
                 is graded if and only if the corresponding Lie
                 algebra homomorphism is graded.

\subsection{The  Baker-Campbell-Hausdorff  formula}\label{BCH}

Let  $G$ be a simply connected nilpotent Lie group with Lie algebra $\mathcal{G}$.
  The exponential map  $\text{exp}:  {\mathcal{G}}\ra G$ is a diffeomorphism.  
   One can then
     pull back the group operation from $G$ to get a group stucture 
 on $\mathcal{G}$.     This group structure can be described by the   Baker-Campbell-Hausdorff formula
 (BCH formula in short),  which expresses the product $X*Y$ ($X, Y\in {\mathcal{G}}$)
    in terms of the iterated Lie brackets  of $X$ and $Y$. 
     The group operation in $G$ will be denoted by $\cdot$.  
   The pull-back group operation   $*$  on $\mathcal{G}$ is defined as follows. 
      For $X,   Y\in \mathcal{G}$,   define
   $$X*Y=\text{exp}^{-1}(\text{exp} X\cdot \text{exp} Y).$$
  Then the first a few terms of the BCH formula (\cite{CG},    page 11)  is given by:
$$ X*Y =X+Y+\frac{1}{2}[X,Y]+\frac{1}{12}[X,[X,Y]]
-\frac{1}{12}[Y, [X, Y]]+\cdots. $$

     Next we define homogeneous distances  on Carnot groups.  With the help of BCH formula, 
  it is often   more  convenient to work with homogeneous distances 
 than with Carnot metrics.  Let   $G$ be  a   Carnot group with  Lie algebra   
$\mathcal
G=V_1\oplus V_2\oplus\cdots \oplus V_r$.  
  Write $x\in \mathcal G$ as $x=x_1+\cdots+ x_r$ with $x_i\in V_i$.  
  Fix a norm   $|\cdot|$ on each layer.  Define a norm $||\cdot||$ on $\mathcal G$
   by:
$$||x||=\sum_{i=1}^r |x_i|^{\frac{1}{i}}.$$
 Now define a homogeneous distance on $G=\mathcal G$  by:  $d(g,h)=||(-g)*h||$. 
  An important fact is that $d$ and $d_c$ are biLipschitz equivalent.   That is, there is a constant $C\ge 1$ such that    $d(p,q)/C\le d_c(p,q)\le C\cdot   d(p,q)$ for all $p, q\in G$. 
  It is often possible to  calculate  or estimate  $d$  by using  the BCH formula.
     Since we are only concerned with quasiconformal maps and biLipschitz maps, 
       it does not matter whether we use $d$ or $d_c$.

\subsection{Quasisymmetric   maps}\label{maps}

Here we recall the definition   of  quasisymmetric map.

Let $\eta: [0,\i)\ra [0,\i)$ be a homeomorphism.
    A homeomorphism
$F:X\to Y$ between two metric spaces is
\e{$\eta$-quasisymmetric} if for all distinct triples $x,y,z\in X$,
we have
\[
   \frac{d(F(x), F(y))}{d(F(x), F(z))}\le \eta\left(\frac{d(x,y)}{d(x,z)}\right).
\]
  If $F: X\ra Y$ is an $\eta$-quasisymmetry, then
  $F^{-1}:    Y\ra X$ is an $\eta_1$-quasisymmetry, where
$\eta_1(t)=(\eta^{-1}(t^{-1}))^{-1}$. See \cite{V}, Theorem 6.3.
 A    homeomorphism between metric spaces
  is quasisymmetric if it is $\eta$-quasisymmetric for some $\eta$. 

We remark that 
     quasisymmetric homeomorphisms between general metric spaces
  are quasiconformal. In the case of Carnot groups  (and more
  generally  Loewner  spaces),  a   homeomorphism  is quasisymmetric if and only
  if it is quasiconformal, see \cite{HK}.

 The main result in \cite{BKR} says that a quasiconformal   map  between two proper, 
locally Ahlfors   $Q$-regular  ($Q>1$)   metric spaces  is absolutely 
 continuous on almost every curve.   This result
   applies to quasiconformal maps    $F:  G\ra G$  on  Carnot groups. 

Pansu \cite{P} proved that a quasisymmetric map 
  $F: G_1\ra G_2$  between two Carnot groups is absolutely continuous:    
  a  measurable set $A\subset G_1$ has measure $0$ if and only if $F(A)$ has measure $0$.

Let $g: (X_1, \rho_1)\ra (X_2, \rho_2)$ be a   homeomorphism   between two
metric spaces.
We define for every $x\in X_1$  and $r>0$,
\begin{align*}
   L_g(x,r)&=\sup\{\rho_2(g(x), g(x')):   \rho_1(x,x')\le r\},\\
   l_g(x,r)&=\inf\{\rho_2(g(x), g(x')):   \rho_1(x,x')\ge r\},
\end{align*}
and set
\[
   L_g(x)=\limsup_{r\ra 0}\frac{L_g(x,r)}{r}, \ \
   l_g(x)=\liminf_{r\ra 0}\frac{l_g(x,r)}{r}.
\]
  Then
\b{equation}\label{e3.001}
  L_{g^{-1}}(g(x))=\frac{1}{l_g(x)} \ \text{ and }\ l_{g^{-1}}(g(x))=\frac{1}{L_g(x)}
\end{equation}
for any $x\in X_1$. If $g$ is an $\eta$-quasisymmetry, then
\b{equation}\label{e3.0}
L_g(x,r)\le \eta(1)l_g(x, r)
\end{equation}
 for all $x\in X_1$ and $r>0$. Hence
if in addition
\[
    \lim_{r\ra 0}\frac{L_g(x,r)}{r}\ \ {\text{or}} \ \ \lim_{r\ra 0}\frac{l_g(x,r)}{r}
\]
exists, then
\[
    0\le l_g(x)\le L_g(x)\le \eta(1) l_g(x)\le \infty.
\]

\subsection{Pansu differentiability  theorem }\label{pansud}

First the definition:

\b{Def}\label{pansu-d}
 Let $G$ and $G'$
  be two Carnot groups endowed with Carnot metrics,  and $U\subset G$, $U'\subset G'$ open subsets.
   A map $F: U\ra U'$ is \e{Pansu  differentiable}
    at $x\in U$  if there exists a graded  homomorphism
     $L: G\ra G'$ such that
     $$\lim_{y\ra x}\frac{d(F(x)^{-1}*F(y),\, L(x^{-1}*y))}{d(x,y)}=0.$$
          In this case, the graded  homomorphism
     $L: G\ra G'$  is called the \e{Pansu  differential} of $F$ at $x$, and
     is denoted by   $dF(x)$.
\end{Def}

We have the following  chain rule for Pansu differentials:

\b{Le} \e{(Lemma  3.7 in \cite{CC})}  \label{chain}
   Suppose $F_1: U_1\ra U_2$ is Pansu differentiable at $p$, and $F_2: U_2\ra U_3$ is 
  Pansu differentiable at $F_1(p)$.  Then $F_2\circ F_1$ is Pansu differentiable at $p$  and 
  $d(F_2\circ F_1)(p)=dF_2(F_1(p))\circ dF_1(p)$. 
\end{Le}

Notice that the Pansu differential of the identity map   $U_1\ra U_1$    is the identity isomorphism.
   Hence if $F: U_1\ra U_2$ is bijective, $F$ is Pansu differentiable at $p\in U_1$  and $F^{-1}$ is Pansu differentiable at $F(p)$,  then 
  $dF^{-1}(F(p))=(dF(p))^{-1}$.

The following  result (except the terminology)   is due to  Pansu
[P].

\b{theorem}\label{pan}
  Let $G, G'$ be   Carnot groups, and $U\subset G$, $U'\subset G'$  open subsets.  
 Let $F: U\ra U'$ be a quasiconformal  map.
   Then $F$ is a.e. Pansu  differentiable. Furthermore, at a.e.
   $x\in   U$, the Pansu  differential $dF(x): G\ra G'$ is a graded
   isomorphism.

\end{theorem}

 In  Theorem \ref{pan}  and the proofs below,  \lq\lq a. e." is with
 respect to the   Lesbegue  measure on $\mathcal{G}=G$.

\section{Quasiconformal  implies    biLipschitz}\label{QCBL}


In this Section, we prove Theorem    \ref{main2} and Theorem \ref{main3}.  They    
 provide   very useful  sufficient conditions for a Carnot group  to 
  be quasisymmetrically rigid. 


We shall need the following result.

\b{Prop}\label{p3} \e{(Proposition   3.4  in \cite{X3}) }
 Let  $G$  and $G'$ be two Carnot groups
  with Lie algebras  $\mathcal G=V_1\oplus \cdots \oplus V_m$  and $\mathcal G'=V'_1\oplus \cdots \oplus  V'_n$ respectively.    Let     $W_1\subset V_1$,  $W'_1\subset V'_1$ be 
  subspaces. Denote by $\mathcal W\subset \mathcal{G}$ and 
$\mathcal {W'}\subset \mathcal{G}'$   respectively the Lie subalgebras
  generated by $W_1$ and $W'_1$.    Let $W\subset G$ and $W'\subset G'$ respectively  be the 
  connected   Lie
  subgroups of $G$ and $G'$    corresponding to 
  $\mathcal W$  and    $\mathcal {W'}$.  
 Let $F: G\ra G'$ be a quasisymmetric  homeomorphism.  
  If     $dF(x)(W_1)\subset W'_1$ for a.e. $x\in G$,   
     then  
   $F$ sends left cosets of $W$ into left  cosets of $W'$.  

  \end{Prop}

Let   $N$ be a Carnot group with Lie algebra  
$\mathcal N=V_1\oplus \cdots \oplus V_r$.  
For any linear subspace $W_1\subset V_1$, let  $\mathcal W$ be the subalgebra of $\mathcal{N}$ generated by $W_1$.  
Denote by $W$ the connected subgroup of $N$   with Lie algebra  $\mathcal W$.  
 Notice that $\mathcal W$ is also a Carnot algebra and can be written as 
   $\mathcal W=W_1\oplus \cdots \oplus  W_r$.    In general, there is some 
  integer  
 $1\le s\le  r$ such that  $W_s\not=0$ and $W_j=0$ for $j> s$.

The following result is very useful for induction argument.

\b{Th}\label{induct}
 Let  $W_1\subset V_1$ be   a non-trivial proper subspace that is invariant under the action of 
 $\text{Aut}_g(\mathcal{N})$. 
Suppose $W$ has the following property:   every 
 $\eta$-quasisymmetric map   $G: W\ra W$ is  a  $(K, C)$-quasi-similarity, where $K$ depends only on $\eta$.     Then $N$  has the same property (possibly with a different constant $K'$ depending only on $\eta$). 


\end{Th}

\b{proof}
The proof is similar to  the arguments in the proofs of 
 Lemma  3.9  and Lemma 3.10 in \cite{X2}. We  include it here 
  mainly   for completeness. 
  Let $F:  N\ra N  $  be an $\eta$-quassymmetric map. 
  By Proposition \ref{p3},  $F$ sends left cosets of $W$ to left cosets of $W$. 
     For each left coset $L$ of $W$,  the restriction $F|_L: L\ra F(L)$ is also an
      $\eta$-quasisymmetric map. 
  Let $L_1=n_1*  W$,   $L_2=n_2*W$ be two left cosets of $W$.     By assumption,  $F|_{L_i}$ is   a  $(K, C_i)$-quasi-similarity, where  $K\ge 1$ and $C_i>0$.   Here    $K$ depends only on $\eta$.  We 
 shall first show that  $C_1/ {2K^2}\le C_2\le 2K^2 C_1$.   
  Suppose $C_2> 2K^2 C_1$.     We will get a contradiction from this. 
 Fix some $w\in W_1$   with $|w|=1$  and consider the image of the  two horizontal lines $n_1* tw $ ($t\in \R$),
  $n_2* tw$ ($t\in \R$)  under $F$.  
   First of all,  by using BCH formula it is easy to show that 
 $d(n_1*tw, n_2* tw)$  is bounded from above by a sublinear function of $t$ as $t\ra \infty$. 
   Since $d(n_1, n_1*tw)=t$, we have 
 $$\frac{d(n_1* tw, n_2*tw)}{d(n_1*tw, n_1)}\ra 0$$
  as $t\ra \infty$.   The quasisymmetric condition   now implies 
 $$\frac{d(F(n_1* tw),   F(n_2*tw))}{d(F(n_1*tw),   F(n_1))}\ra 0. $$
 On the other hand, since $F|_{L_i}$ is a  $(K, C_i)$  quasi-similarity,   we have
  $$d(F(n_2), F(n_2*tw))\ge C_2/K \cdot d(n_2, n_2*tw)>2KC_1\cdot t$$
  and 
  $$d(F(n_1), F(n_1*tw))\le KC_1 \cdot d(n_1, n_1*tw)=KC_1\cdot t.$$  
  It follows that
\b{align*}
\frac{d(F(n_1* tw),   F(n_2*tw))}{d(F(n_1*tw),   F(n_1))} & \ge \frac{d(F(n_2*tw), F(n_2))-d(F(n_2), F(n_1))-d(F(n_1), F(n_1*tw))} { KC_1 t}\\
   & \ge \frac{2KC_1t- d(F(n_2), F(n_1))  -KC_1t}{KC_1t}   \ra 1. 
\end{align*}
  The contradiction shows 
    $C_2\le 2K^2 C_1$.     Similarly,    $C_1/ {2K^2}\le C_2$.

We next show that $F$ is a quasi-similarity.   Fix any left coset $L$ of $W$.  
     Then,  as indicated above,       $F|_L$ is a  $(K, C_0)$-quasi-similarity   for some $C_0>0$.  
  Let $p, q\in N$ be arbitrary points.    Then there are left cosets $L_1$, $L_2$ of $W$ with $p\in L_1$, $q\in L_2$.     By the preceding paragraph, 
  $F|_{L_1}$ is a $(2K^3,   C_0)$-quasi-similarity.  
     Pick   $q'\in L_1$  with $d(p,q')=d(p,q)$.   The quasisymmetry condition implies that 
   $$d(F(p), F(q))\le \eta(1)\cdot   d(F(p), F(q'))\le \eta(1) \cdot 2K^3\cdot C_0 \cdot  d(p, q')=2K^3   \eta(1) \cdot C_0 \cdot d(p, q).$$
    Now the same argument    applied to $F^{-1}$  finishes the proof.

\end{proof}

\b{Le}\label{normal}
 Suppose $W_1\subset V_1$ is a non-trivial proper subspace that is invariant under the action of 
 $\text{Aut}_g(\mathcal{N})$.  If  $\mathcal W=<W_1>$ is  an ideal in $\mathcal{N}$,  then 
 every  $\eta$-quasisymmetric  map $F: N\ra N$ is    a  $(K, C)$-quasi-similarity, where $K$ depends only on $\eta$.  

\end{Le}

\b{proof}
  Since $\mathcal W$ is an ideal in $\mathcal N$,  the subgroup  $W$ corresponding to $\mathcal W$ is a 
connected normal subgroup of $N$.  
By    Proposition \ref{p3},  $F$  maps cosets of $W$ to cosets of $W$.  
Since $W$ is normal in $N$,    any two  cosets 
  $n_1*W$, $n_2 *W$  of $W$ are \lq\lq parallel" in the following sense:
  for any   $p\in n_1 *W$,  we have $d(p,  n_2*W)=d(n_1*W, n_2  *W)$.   Furthermore, the quotient 
 $N/W$ is   a  Carnot group 
  with Lie algebra $V_1/{W_1}\oplus \cdots \oplus  V_r/{W_r}$,   
 and the  distance between $n_1*W $ and $n_2*W $ in $N$ 
 equals the distance between the  two points $n_1W$, $n_2 W$ in $N/W$.  
    In particular,  the metrics on the quotient and each coset   are  geodesic. 
   Now the argument in 
  \cite{SX}    implies that 
 $F$ is   a  $(K, C)$-quasi-similarity, where $K$ depends only on $\eta$.

\end{proof}

\b{Th}\label{normalizer}
Suppose $W_1\subset V_1$ is a non-trivial proper subspace that is invariant under the action of 
 $\text{Aut}_g(\mathcal{N})$.  Let  $N(W_1)=\{v\in V_1:    [v,  W_1]\subset  {\mathcal W}\}$.   
  If $N(W_1)\not= W_1$,   
  then 
 every   $\eta$-quasisymmetric  map $F: N\ra N$ is a   $(K, C)$-quasi-similarity, where $K$ depends only on $\eta$.    
\end{Th}

\b{proof}
   Notice that   $N(W_1)$  is also invariant under the action of 
$\text{Aut}_g(\mathcal{N})$.    Let $N_2\subset N$ be the connected subgroup of $N$ with   Lie algebra 
 $<N(W_1)>$.     
Let $F:  N\ra N$ be an $\eta$-quasisymmetric  map.  
  Proposition \ref{p3}  implies that  $F$ sends left cosets of $N_2$ to left cosets of $N_2$  and  also sends left cosets of $W$ to left cosets of $W$. 
  Clearly,    for each left coset $L$ of  $N_2$,  
 the restriction  $F|_L: L\ra  F(L)$  is also $\eta$-quasisymmetric.   
Notice  that  $\mathcal W$  is  an ideal of $<N(W_1)>$.  
   Since $N(W_1)\not= W_1$,    (the proof of)
 Lemma  \ref{normal}
  yields that     
 the restriction $F|_L$ is a $(K, C)$-quasi-similarity,
  where the constant $K$ depends only on $\eta$.
 Now the  Theorem  follows from     Theorem   \ref{induct}.

\end{proof}

\section{Invariant subspaces}\label{invariant}

In this Section we first explain how to find invariant subspaces, 
  then  introduce some constructions that produce Carnot algebras with 
  reducible first layer.


\subsection{Rank  and invariant subspaces}\label{rank}

In this  Subsection we  define the  rank of elements in a Lie algebra,  
   and  explain how to use rank to  find invariant subspaces.   The notion of rank    appeared in \cite{O}.     There is a  characterization of non-rigid Carnot groups in terms of rank \cite{OW}.

Let  $f:   X\ra  Y$  be  a     linear   map between two vector spaces.   The 
  {\it rank}   $\text{rank}(f)$   of
  $f$ is defined to be the dimension of the image $f(X)$  of $f$.   
  Recall   that    for any positive integer $m$,   $f$ induces a  linear map  $\wedge^{m}f:   
 \wedge^{m}  X \ra \wedge^m  Y$  between the $m$-fold exterior products.

\b{Le}\label{rankc}
Let $k$ be a non-negative integer and $f:   X \ra   Y $ be a  linear map. 
  Then $\text{rank}(f)\le k$ if and only if $\wedge^{k+1}f=0$.

\end{Le}

\b{proof}
  First assume $\text{rank}(f)\le k$.  Then $f(X)$ has dimension at most $k$.  
  Hence for any $x_1, \cdots, x_{k+1}\in  X$,  
 we have 
 $$\wedge^{k+1}f(x_1\wedge \cdots \wedge x_{k+1})=f(x_1)\wedge \cdots\wedge f(x_{k+1})=0.$$

Conversely,   assume $\text{rank}(f)\ge k+1$.  Then there are vectors $x_1,  \cdots,  x_{k+1}\in X$ such that    $f(x_1), \cdots,  f(x_{k+1})$  are linearly independent.  
  Hence 
$$\wedge^{k+1}f(x_1\wedge \cdots \wedge x_{k+1})=f(x_1)\wedge \cdots\wedge f(x_{k+1})\not=0.$$

\end{proof}

For an element $x\in \mathcal{N}$ in a  Lie algebra, let $\text{rank}(x)$ be the rank of 
 the linear transformation  $ad(x):  \mathcal{N}\ra   \mathcal{N}$,    $ad(x)(y)=[x,y]$.   In other words, 
  $\text{rank}(x)$  is the dimension of the image of $ad(x)$. 
   For any isomorphism $A: \mathcal N_1 \ra \mathcal N_2$ of  Lie algebras,  
  we have $\text{rank}(x)=\text{rank}(A(x))$  for any $x\in \mathcal N_1$.  
The rank of a Lie algebra  $\mathcal N$
 is   defined  by 
   $$\text{rank}(\mathcal N)=\min\{\text{rank}(x):  0\not=x\in \mathcal N\}. $$
  For any $\lambda\not=0$  and any $x\in \mathcal{N}$,   we clearly have 
 $\text{rank}(\lambda x)=\text{rank}(x)$.  Hence the rank induces a function 
   $R:   P\mathcal{N}  \ra  \Z$    on the projective space
  $P\mathcal{N}$:
  $$R([x])=\text{rank}(x)\;  \text{for}\;   0\not=x\in \mathcal{N}.$$ 
Set  $A_k=\{[x]\in P\mathcal{N}:   \text{rank}(x)\le k\}$.   Clearly $A_k\subset A_{k+1}$ for all $k$.

  Let $\{e_1, \cdots, e_n\}$ be a   basis of   a  Lie algebra $\mathcal{N}$.  Then any element 
$x\in \mathcal{N}$  can be written as:      $x=x_1 e_1+ \cdots  +x_n e_n$  ($x_i\in \R$).
   Then $ad(x)=  x_1 ad(e_1)+\cdots   +x_n ad(e_n)$.  It follows that 
 $$\wedge^{m} ad (x)=\sum_{i_1, \cdots, i_m} x_{i_1}\cdots x_{i_m} ad(e_{i_1})\wedge \cdots  \wedge  ad(e_{i_m})$$
  is a linear combination of the linear maps $ ad(e_{i_1})\wedge \cdots \wedge  ad(e_{i_m})$
  and the coefficients are polynomials in $x_1,  \cdots,  x_n$.  
  Now it is easy to see by using Lemma \ref{rankc}  that $A_k\subset   P\mathcal{N}$  is  an  algebraic    variety.

Now let $\mathcal N=V_1\oplus \cdots \oplus V_r$ be a Carnot algebra. 
   Define
  $$r_1(\mathcal N)=\min\{\text{rank}(x):  0\not=x\in  V_1\}. $$
  The quantity $r_1(\mathcal N)$ will be called the {\it rank} of $\mathcal N$ and of the 
    corresponding Carnot group $N$. 
  Let $W_1\subset  V_1$ be the subspace of $V_1$ spanned by elements 
 $x\in V_1$ satisfying $\text{rank}(x)=r_1(\mathcal N)$.  
  Since graded isomorphisms preserve  the first layer $V_1$,  we see that $W_1$ is a 
 non-trivial subspace of $V_1$ invariant under the action of 
 $\text{Aut}_g(\mathcal N)$.    Very often $W_1$ is a proper subspace of $V_1$. For example,   this is the case for most non-rigid Carnot groups \cite{X1}.
    Similarly,  for $1\le i\le r-1$, define
 $$r_{1,i}(\mathcal N)=\min\{\text{rank}(ad (x)|_{V_i}):  0\not=x\in  V_1\}. $$
Let $W_{1,i}\subset  V_1$ be the subspace of $V_1$ spanned by elements 
 $x\in V_1$ satisfying $\text{rank}(ad (x)|_{V_i})=r_{1,i}(\mathcal N)$.  
   Then  $W_{1,i}$ is  also   a 
 non-trivial subspace of $V_1$ invariant under the action of 
 $\text{Aut}_g(\mathcal N)$.

Since $\{A_k\}$ is a nested sequence of algebraic varieties, the best 
  chance to obtain proper invariant subspace is to consider elements of minimal rank.
   This is why we defined the subspaces $W_1$ and $W_{1, i}$.  

Recall that a Carnot group is called rigid if the space of smooth   contact maps is finite dimensional.  
   Ottazzi  (\cite{O})  showed that   if  $N$  has rank  at most one, then $N$ is non-rigid.  
 Furthermore,  Ottazzi  and Warhurst  (\cite{OW})  proved that a Carnot group is non-rigid 
 if and only if  the complexification     $\mathcal N\otimes  \C$ of its Lie algebra
 has rank at most one (as a complex Lie algebra).

\subsection{Product constructions}\label{product}

In this Subsection we introduce a few product  type   constructions for Carnot algebras and Carnot groups. They produce Carnot algebras with proper invariant subspaces.

{\bf{Direct Product}}. 
The  basic product construction is the direct product of two Lie algebras. We will use the direct sum notation.
   Let $\mathcal{N}=V_1\oplus \cdots \oplus V_m$ and $\mathcal{N}'=V'_1\oplus \cdots  \oplus  V'_n$  be two     Carnot algebras. We may assume $m\ge n$.   We will write
  $\mathcal{N}'=V'_1\oplus \cdots  \oplus  V'_m$  with $V'_i=\{0\}$ for $i> n$. 
  Then   $\mathcal{N}\oplus \mathcal{N}'=(V_1\oplus V'_1)\oplus \cdots \oplus (V_m\oplus V'_m)$
  is also a  Carnot algebra.  
  Notice that if $x_1\in V_1$ and $x'_1\in   V'_1$,   then 
   $\text{rank}(x_1+x'_1)=\text{rank}(x_1)+  \text{rank}(x'_1)$.  It follows that,
  if $r_1(\mathcal N)\not=r_1(\mathcal N')$,  then the subspace of $V_1\oplus V'_1$ spanned by elements of minimal rank  is a non-trivial proper subspace   invariant 
 under the action of  $\text{Aut}_g(\mathcal{N}\oplus \mathcal{N}')$.

{\bf{Central  Product}}. 
For any ideal  $\mathcal{I}\subset  \mathcal{N}\oplus \mathcal{N}'$,  the quotient 
  $ (\mathcal{N}\oplus \mathcal{N}')/\mathcal{I}$  is a  Lie algebra.   In general, this quotient is not a Carnot algebra.    An ideal     $\mathcal{I}\subset  \mathcal{N}\oplus \mathcal{N}'$  is called a  graded ideal if 
     $\mathcal{I}=[\mathcal{I}\cap (V_1\oplus V'_1)]\oplus  \cdots  \oplus [\mathcal{I}\cap (V_m\oplus V'_m)]$.
   If   $\mathcal{I}\subset  \mathcal{N}\oplus \mathcal{N}'$  is  a  graded ideal,  then  
   the quotient   $ (\mathcal{N}\oplus \mathcal{N}')/\mathcal{I}$     is  also  a Carnot algebra:
      $ (\mathcal{N}\oplus \mathcal{N}')/\mathcal{I}=V''_1\oplus \cdots \oplus V''_m$,
  where 
  $$V''_i=(V_i\oplus V'_i)/   [(V_i\oplus V'_i)\cap  \mathcal{I}].$$

  Let  $\mathcal{I}\subset V_m\oplus V'_m$   be  any   linear subspace.
   Since  $V_m\oplus V'_m$  lies in the center of $  \mathcal{N}\oplus \mathcal{N}'$, 
    we see that 
$\mathcal{I}$  is     a graded ideal 
  of $  \mathcal{N}\oplus \mathcal{N}'$.    So the associated quotient is a  Carnot algebra. 
A  particular case is the   following.  Let $W\subset V_m$,  $W'\subset V'_m$ be linear subspaces  and 
 $f:  W\ra W'$  a  linear isomorphism.    Then 
 $\mathcal{I}=\{f(w)-w|   w\in W\}$  is a   graded ideal of  $  \mathcal{N}\oplus \mathcal{N}'$.
  The associated quotient is called a {\it central product} of the two   Carnot  algebras
  $\mathcal{N}$  and $\mathcal{N}'$.

Let  $\ol {\mathcal N}$ be a central product of $\mathcal N$ and $\mathcal N'$.
   Notice that   the compositions
$\mathcal N\hookrightarrow \mathcal N\oplus \mathcal N'\ra \ol {\mathcal N}$,
${\mathcal N'}\hookrightarrow \mathcal N\oplus \mathcal N'\ra \ol {\mathcal N}$
   are  injective,  and $[V_1, \mathcal N']=0$, $[V'_1, \mathcal N]=0$.
    It follows that each  $x_1\in V_1$ has the same rank  as an element in $\mathcal N$ and an element in $\ol {\mathcal N}$.    Similarly  for $x'_1\in V'_1$.   Furthermore,
  if $x_1\in V_1$ and $x'_1\in V'_1$,   then 
   $\text{rank}(x_1+x'_1)\ge \max\{ \text{rank}(x_1),   \text{rank}(x'_1)\}$.
   Consequently,  
 if $r_1(\mathcal N)\not=r_1(\mathcal N')$,  then the subspace of $V_1\oplus V'_1\subset \ol {\mathcal N}$ spanned by elements of minimal rank  is a non-trivial proper subspace   invariant 
 under the action of  $\text{Aut}_g(\ol{\mathcal{N}})$.
 We have observed that both $\mathcal  N$ and $\mathcal N'$ are 
isomorphic to their images 
 in the   central  product. 

Given an  $m$-step Carnot algebra $\mathcal N$ with $m\ge 2$,  
  one can   always find another $m$-step Carnot algebra $\mathcal N'$ with 
 $r_1(\mathcal N')>r_1(\mathcal N)$. For instance, one may choose an
  $m$-step free nilpotent  Lie algebra on sufficiently many generators.  
Hence, by the above discussion, $\mathcal N$ always embeds in a Carnot algebra 
  with a reducible first layer. 

{\bf{Level  One    Product}}. 
Now let $\mathcal N=V_1\oplus V_2$ and $\mathcal N'=V'_1\oplus  V'_2$ be $2$-step Carnot algebras.  Let $X_0\in V_1\backslash\{0\}$ and 
$X'_0\in V'_1\backslash\{0\}$.   Let     $U_1\subset V_1$ 
  be  a   subspace 
complementary to 
 $\R X_0$  and
      $U'_1\subset V'_1$   a subspace  complementary to $\R X'_0$. 
  We define a new $2$-step Carnot algebra 
 $\ol {\mathcal N}={\ol V_1}\oplus {\ol V_2}$ as follows.
    The first layer ${\ol V_1}=U_1\oplus U'_1\oplus \R X$, where  $\R X$ is a $1$-dimensional vector space with basis vector $X$. 
   The second layer ${\ol V_2}=V_2\oplus V'_2$. 
  The bracket on $\ol {\mathcal N}$ is defined as follows.  
 For any $u_1, u_2\in U_1$, $u'_1, u'_2\in U'_1$, $a, b\in \R$,
 define
 \b{align*}
& [u_1+u'_1+a X, \;  u_2+u'_2+b X]\\
= &[u_1, u_2]+b\cdot [u_1,  X_0]+[u'_1, u'_2]+b\cdot [u'_1, X'_0]+a\cdot [X_0,   u_2] +a\cdot [X'_0,   u'_2].
\end{align*}
  We call  $\ol {\mathcal N}$  a 
   {\it  level one product}    of $\mathcal N$ and $\mathcal N'$. 
  It depends on  the choices of  $X_0$, $X'_0$, $U_1$, $U'_1$.   
We notice that both $\mathcal N$ and $\mathcal N'$ are isomorphic to their images 
 in the level one product. 

For the following   Lemma,   notice that we  can always choose $X_0$, $X'_0$ and $U_1$, $U'_1$ so that $U_1$ and $U'_1$  contain elements of minimal rank.

\b{Le}\label{level1}
Let $\ol{\mathcal N}$ be  a level one product of $\mathcal N$ and $\mathcal N'$ defind as above.    Suppose there are $\tilde u_1\in U_1$ and $\tilde u'_1\in U'_1$ such that
   $\text{rank}(\tilde u_1)=r_1(\mathcal N)$ and 
 $\text{rank}(\tilde u'_1)=r_1(\mathcal N')$.
  If   $r_1(\mathcal N)\ge 1$  and    $r_1(\mathcal N')\ge 1$,  
    then    $\ol{\mathcal N}$  has a reducible first layer. 
\end{Le}

\b{proof}
To avoid confusion, if $v\in \ol {V_1}$,  we use $\ol{r}(v)$ to denote the rank of 
 $v$ in the Lie algebra $\ol{\mathcal N}$.  We still use $\text{rank}(x)$ to denote rank of elements  $x$  in the Lie algebras $\mathcal N$ and $\mathcal N'$. 
The definition of the Lie bracket in  $\ol{\mathcal N}$  shows that
 for any $u_1\in U_1$, $u'_1\in U'_1$, we have
 $[u_1,   \ol{V_1}]=[u_1, V_1]$ and  $[u'_1,   \ol{V_1}]=[u'_1, V'_1]$.  This implies
  $\ol{r}(u_1)=\text{rank}(u_1)$ and 
$\ol{r}(u'_1)=\text{rank}(u'_1)$.      By the  assumption, there are $\tilde u_1\in U_1$ and $\tilde u'_1\in U'_1$ such that
   $\ol{r}(\tilde u_1)=r_1(\mathcal N)$ and 
 $\ol{r}(\tilde u'_1)=r_1(\mathcal N')$.
  We next show that if $v\in \ol{V_1}\backslash (U_1\oplus U'_1)$, then 
  $\ol{r}(v)\ge r_1(\mathcal N) +r_1(\mathcal N')$.  
Since we assume  $r_1(\mathcal N)\ge 1$   and    $r_1(\mathcal N')\ge 1$,  
  this implies that  all elements of minimal rank are contained in 
 $U_1\oplus U'_1$, and the Lemma follows. 

   Let $v\in \ol{V_1}\backslash (U_1\oplus U'_1)$.  Write $v$ as 
  $v=u_1+u'_1+a X$ with $u_1\in U_1$, $u'_1\in U'_1$ and $a\not=0$.  
    We may assume $a=1$. So $v=u_1+u'_1+ X$.  
     Then 
 $$[v, \ol{V_1}]=\big\{[u_1+X_0, u_2]+b [u_1,   X_0]+[u'_1+X'_0,  u'_2]+b[u'_1,  X'_0]:  u_2\in U_1,  u'_2\in U'_1, b\in \R\big\}.$$
 Let $P:  \ol{V_2}=V_2\oplus V'_2\ra V_2$  be the projection.  
 Then 
 \b{align*}
P([v, \ol{V_1}])& =\{[u_1+X_0, u_2]+b [u_1,   X_0]:   u_2\in U_1,   b\in \R\}\\
  &=\{[u_1+X_0, u_2]+ [u_1+X_0,  b X_0]:   u_2\in U_1,   b\in \R\}\\
&=\{[u_1+X_0, u_2+b X_0]:   u_2\in U_1,   b\in \R\}\\
&=[u_1+X_0, V_1] .
\end{align*}
  Hence $\dim(P([v, \ol{V_1}]))= \text{rank}(u_1+X_0)\ge r_1(\mathcal N)$.  
  On the other hand,
   \b{align*}
 &[v, \ol{V_1}]\cap P^{-1}(0)\\
  =&\{[u'_1+X'_0,  u'_2]+b[u'_1,  X'_0]: 
[u_1+X_0, u_2]+b [u_1,   X_0]=0,   u_2\in U_1,  u'_2\in U'_1, b\in \R\}\\
=&\{[u'_1+X'_0,  u'_2+bX'_0]: 
[u_1+X_0, u_2]+b [u_1,   X_0]=0,   u_2\in U_1,  u'_2\in U'_1, b\in \R\}.
  \end{align*}
  Notice that there is no restriction on $u'_2$. 
  Also, for any $b\in \R$,  the  equality 
 $[u_1+X_0, u_2]+b [u_1,   X_0]=0$ holds for $u_2=b u_1$.  
 Hence 
 $[v, \ol{V_1}]\cap P^{-1}(0)=[u'_1+X'_0,   V'_1]$.
  It follows that  
$$\dim([v, \ol{V_1}]\cap P^{-1}(0))= \text{rank}(u'_1+X'_0)\ge  r_1(\mathcal N')$$ 
  and so 
 $$\ol{r}(v)=\dim([v, \ol{V_1}]) =\dim(P([v, \ol{V_1}]))  +
\dim([v, \ol{V_1}]\cap P^{-1}(0))
\ge   r_1(\mathcal N)  +r_1(\mathcal N').$$


\end{proof}

\section{QC maps on $2$-step Carnot groups}\label{2stepcase}

In this Section we prove Theorem \ref{main1}.



{\bf Throughout this Section,   $N$ will be a $2$-step Carnot group with Lie algebra 
$\mathcal{N}=V_1\oplus V_2$,   $W_1\subset V_1$  a 
non-trivial proper subspace  
   invariant under  the action of $\text{Aut}_g(\mathcal{N})$  on $V_1$,
  and $F: N\ra N$  an $\eta$-quasisymmetric map. }

   Let  $\mathcal W$ be the subalgebra of $\mathcal{N}$ generated by $W_1$,
  and   $W$ the connected subgroup of $N$   with Lie algebra  $\mathcal W$.  
 Notice that $\mathcal W$ is also a Carnot algebra and can be written as 
   $\mathcal W=W_1\oplus  W_2$.     We may have $W_2=0$.  This happens if and only if $W_1$ is abelian;   that is,  $[W_1, W_1]=0$.
Let $\tilde W_1\subset V_1$ be a   subspace complementary   to  $W_1$,    and $\tilde W_2\subset V_2$ be a subspace complementary to  $W_2$.  Then we have $V_1=W_1\oplus \tilde W_1$ and $V_2=W_2\oplus \tilde W_2$.    We denote by $P_2: V_2\ra W_2$ and  $\tilde P_2: V_2\ra \tilde W_2$ respectively   the projections with respect to the 
 direct sum decomposition $V_2=W_2\oplus \tilde W_2$.  That is, if 
  $w=w_2+\tilde w_2$ with $w_2\in W_2$ and $\tilde w_2\in \tilde W_2$, then
 $P_2(w)=w_2$ and $\tilde P_2(w)=\tilde w_2$.

Theorem \ref{main1} follows from   Lemma  \ref{normal} in the case when   $W_2=V_2$.  From now on we shall assume 
 $W_2\not=V_2$.   
We fix  inner products on $V_1$ and $V_2$ such that $W_i$ and $\tilde W_i$ are perpendicular.     Then there exists some constant $A\ge 1$ such that 
 $|[v, w]|\le A \cdot |v|\cdot |w|$   for all $v, w\in V_1$.   Fix some $e\in \tilde W_2$ with $|e|=1$.

We recall that,  if $F: X\ra Y$ is $\eta$-quasisymmetric, then $F^{-1}: Y\ra X$ is 
  $\eta_1$-quasisymmetric  with $\eta_1(t)=(\eta^{-1}({t}^{-1}))^{-1}$. 
  Without loss of generality we may assume $\eta(1)\ge 1$. It follows that  $\eta_1(1)\ge 1$.  
 These inequalities will be implicitly used  throughout the paper.


\b{Le}\label{key}
  Suppose $W_2\not=V_2$.   
Let $L$  be   a left coset of $W$.  
   Suppose $p,q\in L$ are such that $l_F(p)>  C_1\cdot  L_F(q)$  with 
$C_1= 200A \eta_1(1)$. 
   Write $q=p*(u_1+u_2)$ with $u_1\in W_1$ and $u_2\in W_2$.   Then for any $s=p*(tu_1+w'_2)$  with $|t|\ge 1$ and $w'_2\in W_2$,    we   have 
 $$ \text{L}_F(s)\le  \frac{2(\eta_1(1))^2} {\sqrt{|t|}}\cdot  \text{L}_F(q).$$

\end{Le}

\b{proof}  Set   $p'=F(p)$,   $q'=F(q)$   and  $L'=F(L)$.  
  The assumption   and  (\ref{e3.001})  imply        $l_{F^{-1}}(q')>  C_1\cdot   L_{F^{-1}}(p')$. 
Let $\{r_j\}$ be an arbitrary sequence of positive reals such that $r_j\ra 0$.
  Then 
\b{equation}\label{e5.10}
\liminf\frac{l_{F^{-1}}(q', r_j)}{r_j}\ge C_1 \cdot \limsup\frac{L_{F^{-1}}(p', r_j)}{r_j}.
\end{equation}
  Fix some  $e\in \tilde W_2$  with $|e|=1$.  
    We shall look at the image of $r^2_j e*L'$ under $F^{-1}$.

Denote $\bar p_j=r^2_j e*p'$ and $\bar q_j=r^2_j e*q'$. 
   Since $d(p', \bar p_j)=d(q',  \bar q_j)=r_j$,  we have 
  $$\frac{d(F^{-1}(\bar q_j), q)}{r_j}\ge \frac{l_{F^{-1}}(q',   r_j)}{r_j}$$
  and 
$$\frac{d(F^{-1}(\bar p_j), p)}{r_j}\le \frac{L_{F^{-1}}(p',   r_j)}{r_j}.$$
  Let $p_j,  q_j$   
be points on $ F^{-1}(r^2_j e*L')$ nearest to $p$ and $q$ respectively. 
  Since $W_2$ and $\tilde W_2$ are perpendicular,   we have $d(p', \bar p_j)=d(p', r^2_j e*L')=r_j$. 
  In particular,  
    $d(p',  \bar p_j)\le d(p',  F(p_j))$.  Now the  quasisymmetry condition of $F^{-1}$  implies  
   $d(p,  F^{-1}(\bar p_j))\le \eta_1(1) d(p,  p_j)$.  Similarly we have
 $d(q,  F^{-1}(\bar q_j))\le \eta_1(1) d(q,  q_j)$.  
 It follows that 
\b{equation}\label{e3.01}
\frac{d(q_j, q)}{r_j}\ge \frac{1}{\eta_1(1)}\cdot \frac{d(F^{-1}(\bar q_j), q)}{r_j}\ge  
\frac{1}{\eta_1(1)}\cdot  \frac{l_{F^{-1}}(q',   r_j)}{r_j}.
\end{equation}
    The choice of $p_j$ implies 
\b{equation}\label{e5.13}
\frac{d(p_j, p)}{r_j}\le 
\frac{d(F^{-1}(\bar p_j), p)}{r_j}\le \frac{L_{F^{-1}}(p',   r_j)}{r_j}.
\end{equation}
  It follows  from (\ref{e3.01}),  
 (\ref{e5.10})   and  (\ref{e5.13})   that 
  \b{align*}
\liminf \frac{d(q_j, q)}{r_j}\ge \liminf   \frac{1}{\eta_1(1)}\cdot  \frac{l_{F^{-1}}(q',   r_j)}{r_j}
  &  \ge    \frac{C_1}{\eta_1(1)}\cdot  \limsup \frac{L_{F^{-1}}(p',   r_j)}{r_j} \\
 & \ge    \frac{C_1}{\eta_1(1)}\cdot   \limsup 
    \frac{d(p_j, p)}{r_j}.
\end{align*}
  Now the    choice of $C_1$   implies that 
\b{equation}\label{e3.1}
   {d(p, p_j)}\le \frac{1}{101A}d(q, q_j)
\end{equation}
  for all sufficiently large $j$.

Next we shall look at $d(p, p_j)$  and $d(q, q_j)$.

Since $q, p$ lie on the same left  coset,  we  can write
  $q=p*(u_1+u_2)$   for some $u_1\in W_1$, $u_2\in W_2$.   
  Similarly,   $q_j=p_j*(w_1+w_2)$  for some $w_1\in W_1$, $w_2\in W_2$.  
  Let $o_j=p_j*(w'_1+w_2')$ ($w_1'\in W_1,  w'_2\in W_2$)  be an arbitrary point on the left coset $p_j*W$. 
  Also 
   write $p_j=p*(x_1+ \tilde{x}_1+x_2+ \tilde{x}_2)$
  with  $x_i\in W_i$,  $\tilde{x}_i\in \tilde W_i$.   Although the  $w_i$'s,   $x_i$'s and
     $\tilde{x}_i$'s      depend on 
 $r_j$,  we shall surpress the dependence to simplify the notation.

Next we calculate $d(p, p_j)$  and $d(q, o_j)$.
  Notice 
$$(-p)*p_j= (-p)*p*(x_1+ \tilde{x}_1+x_2+ \tilde{x}_2)=x_1+ \tilde{x}_1+x_2+ \tilde{x}_2.$$  
So  
  \b{equation}\label{e3.2}
d(p, p_j)=d(o,  (-p)*p_j)=|x_1|+|\tilde x_1|+|x_2|^{\frac{1}{2}}+|\tilde x_2|^{\frac{1}{2}}.\end{equation}
 Using the BCH formula, we obtain: 
\b{align*}
  &(-q)*o_j  \\
& =(-q)*p_j*(w'_1+w'_2)\\
&=(-u_1-u_2)*(-p)*p*(x_1+ \tilde{x}_1+x_2+ \tilde{x}_2)*(w'_1+w'_2)\\
 &=(-u_1-u_2)*(x_1+ \tilde{x}_1+x_2+ \tilde{x}_2)*(w'_1+w'_2)\\
& =(x_1-u_1+w'_1)+\tilde x_1+x_2-u_2+w'_2+\tilde x_2-\frac{1}{2}[u_1, x_1]-\frac{1}{2}[u_1, \tilde x_1]+\frac{1}{2}[x_1-u_1, w'_1]+\frac{1}{2}[\tilde x_1, w'_1]\\
&=(x_1-u_1+w'_1)+\tilde x_1+P_2((-q)*o_j )+\tilde P_2((-q)*o_j),
\end{align*}
   where 
 \b{equation}\label{e3.3}
P_2((-q)*o_j )=x_2-u_2+w'_2-\frac{1}{2}[u_1,  x_1]+\frac{1}{2}[x_1-u_1,  w'_1]+
\frac{1}{2}P_2([\tilde x_1, u_1+w'_1])
\end{equation}
  and  
\b{equation}\label{e3.4}
\tilde P_2((-q)*o_j)=\tilde x_2+\frac{1}{2}\tilde P_2([\tilde x_1, u_1+w'_1]).
\end{equation}
  Hence
 \b{equation}\label{e3.5}
d(q, o_j)=|x_1-u_1+w'_1|+|\tilde x_1|+|P_2((-q)*o_j )|^{\frac{1}{2}}
+|\tilde P_2((-q)*o_j)|^{\frac{1}{2}}.
\end{equation}

Since $q_j$ is a point on $p_j*W$ nearest to $q$, we have 
   $d(q_j, q)\le d(o_j, q)$ for any $o_j\in p_j*W$.   By (\ref{e3.1}), we have 
 $d(p, p_j)\le d(q, o_j)/(101A)$.  
  In particular, this inequality holds for the point $o_j=p_j*(w'_1+w_2')$  where $w'_1=u_1$ and $w'_2$ is such that $P_2((-q)*o_j )=0$.  
  Now  using (\ref{e3.2}), (\ref{e3.4}) and (\ref{e3.5})  we obtain:
\b{align*}
 |x_1|+|\tilde x_1|+|x_2|^{\frac{1}{2}}+|\tilde x_2|^{\frac{1}{2}}& \le \frac{1}{101A}
 \left(|x_1|+|\tilde x_1|+\left|\tilde x_2+\tilde P_2([\tilde x_1, u_1])\right|^{\frac{1}{2}}\right)\\
& \le \frac{1}{101A}
 \left(|x_1|+|\tilde x_1|+|\tilde x_2|^{\frac{1}{2}}+\left|\tilde P_2([\tilde x_1, u_1])\right|^{\frac{1}{2}}\right),
\end{align*}
  which implies
\b{equation}\label{e3.6}
 |x_1|+|\tilde x_1|+|x_2|^{\frac{1}{2}}+|\tilde x_2|^{\frac{1}{2}}\le \frac{1}{100A}
\left|\tilde P_2([\tilde x_1, u_1])\right|^{\frac{1}{2}}.
\end{equation}

Next we consider a point $s\in p*W$ of the form $s=p*(tu_1+u'_2)$, where  $u'_2\in W_2$ is arbitrary and $|t|\ge 1$. 
  We claim that 
\b{equation}\label{e3.6.5}
d(s,\,  p_j*W)\ge  \frac{\sqrt{|t|}}{2} d(q, q_j).
\end{equation} 
In the formulas (\ref{e3.3}),  (\ref{e3.4}) and (\ref{e3.5})  we replace $u_1$ with $tu_1$, 
 $u_2$ with $u'_2$ and $q$ with $s$   to  obtain:
\b{equation}\label{e3.7}
P_2((-s)*o_j )=x_2-u'_2+w'_2-\frac{1}{2}t[u_1,  x_1]+\frac{1}{2}[x_1-tu_1,  w'_1]+
\frac{1}{2}P_2([\tilde x_1, tu_1+w'_1]),
\end{equation}
\b{equation}\label{e3.8}
\tilde P_2((-s)*o_j)=\tilde x_2+\frac{1}{2}\tilde P_2([\tilde x_1, tu_1+w'_1])
\end{equation}
    and
 \b{equation}\label{e3.9}
d(s, o_j)=|x_1-tu_1+w'_1|+|\tilde x_1|+|P_2((-s)*o_j )|^{\frac{1}{2}}
+|\tilde P_2((-s)*o_j)|^{\frac{1}{2}}.
\end{equation}

Now for any $o_j=p_j*(w'_1+w'_2)\in p_j*W$, either $|w'_1-tu_1|\ge \sqrt{|t|}\left|\tilde P_2([\tilde x_1, u_1])\right|^{\frac{1}{2}}$  or  
 $|w'_1-tu_1|\le \sqrt{|t|}\left|\tilde P_2([\tilde x_1, u_1])\right|^{\frac{1}{2}}$.
  If    $|w'_1-tu_1|\ge \sqrt{|t|}\left|\tilde P_2([\tilde x_1, u_1])\right|^{\frac{1}{2}}$,
  then   (\ref{e3.6})   and (\ref{e3.9})  imply
 $$  d(s, o_j)\ge |x_1-tu_1+w'_1|\ge (\sqrt{|t|}-\frac{1}{100A}) \left|\tilde P_2([\tilde x_1, u_1])\right|^{\frac{1}{2}}\ge \frac{\sqrt{|t|}}{2} d(q, q_j).   $$
  If  
$|w'_1-tu_1|\le \sqrt{|t|}\left|\tilde P_2([\tilde x_1, u_1])\right|^{\frac{1}{2}}$,
  then  by  (\ref{e3.6})
\b{align}
|[\tilde x_1,   w'_1-t u_1]|  &  \le A\cdot  |\tilde x_1|\cdot |w'_1-t u_1|\nonumber \\
  &  \le A\cdot  |\tilde x_1|\cdot \sqrt{|t|}\left|\tilde P_2([\tilde x_1, u_1])\right|^{\frac{1}{2}}  \nonumber  \\
  &  \le \frac{1}{100} \sqrt{|t|}\left|\tilde P_2([\tilde x_1, u_1])\right|\label{e3.10}
\end{align}
for all sufficiently large $j$. 
In this case, by 
       (\ref{e3.9}),   (\ref{e3.8})   and then  (\ref{e3.6}),   (\ref{e3.10})
 \b{align*}
 d(s, o_j)  & \ge  |\tilde P_2((-s)*o_j)|^{\frac{1}{2}} \\
&  =|\tilde x_2+\frac{1}{2}\tilde P_2([\tilde x_1, tu_1+w'_1])|^{\frac{1}{2}} \\
  &= |\tilde x_2+\frac{1}{2}\tilde P_2([\tilde x_1, w'_1-t u_1])+  t\tilde P_2([\tilde x_1,  u_1]) |^{\frac{1}{2}} \\
& \ge  \left |\left(|t|-\frac{1}{(100A)^2}-\frac{1}{2}\cdot \frac{1}{100} \sqrt{|t|}\right)  \tilde P_2([\tilde x_1,  u_1]) \right|^{\frac{1}{2}} \\
  & \ge   \frac{\sqrt{|t|}}{2} d(q, q_j)
\end{align*}
  as $|t|\ge 1$.   Hence (\ref{e3.6.5}) holds.

Set $s'=F(s)$.  Fix a sequence $r_j\ra +0$ such that
 $l_{F^{-1}}(s')=\lim_j \frac{l_{F^{-1}}(s',  r_j)}{r_j}$.
  By (\ref{e3.0}),  (\ref{e3.6.5})     and  (\ref{e3.01}) 
\b{align*}
l_{F^{-1}}(s')& =\lim_j \frac{l_{F^{-1}}(s',  r_j)}{r_j}\\
  & \ge \frac{1}{\eta_1(1)}\cdot \limsup_j \frac{L_{F^{-1}}(s', r_j)}{r_j}\\
& \ge \frac{1}{\eta_1(1)}\cdot \limsup_j \frac{d(s,  F^{-1}(r_j^2 e*s'))}{r_j}\\
& \ge \frac{\sqrt{|t|}}{2\eta_1(1)}\cdot \limsup_j \frac{d(q, q_j)}{r_j}\\
& \ge \frac{\sqrt{|t|}}{2(\eta_1(1))^2}\cdot \limsup_j \frac{l_{F^{-1}}(q', r_j)}{r_j}\\
& \ge \frac{\sqrt{|t|}}{2(\eta_1(1))^2}\cdot   l_{F^{-1}}(q').
\end{align*}
   The Lemma now follows from (\ref{e3.001}).

\end{proof}

Let  $\pi_1:  \mathcal N=V_1\oplus V_2 \ra V_1$  be the projection onto  $V_1$. 
 Notice that for each left coset $L=p*W$ of $W$,   the projection  $\pi_1(L)=\pi_1(p)+W_1$    is an affine subspace of $V_1$ parallel to   $W_1$. 
  By \lq\lq a hyperplane $H$  in  $ \pi_1(L)$" we mean an affine subspace $H\subset \pi_1(L)$ such that   $\dim(H)=\dim(\pi_1(L))-1$.

\b{Le}\label{add1}
   For any two points $p, q\in L$   satisfying  
$\pi_1(p)=\pi_1(q)$  we  have   $l_F(p)\le    C_1\cdot  L_F(q)$, where 
$C_1= 200A \eta_1(1)$.

\end{Le}

\b{proof}  We use the notation in the proof of  Lemma \ref{key}.  
   Assume $l_F(p)>    C_1\cdot  L_F(q)$.   Write $q=p*(u_1+u_2)$.  
  Then   equality (\ref{e3.6})  holds,   which  implies $u_1\not=0$.  
  Hence   $\pi_1(p)\not=\pi_1(q)$.  
 

\end{proof}

\b{Le}\label{halfspace}
 Let  $L$ be a left coset of $W$ and $p, q\in L$.  
Suppose $\text{l}_F(p)> (C_2)^{2m} \text{L}_F(q)$,  where $m=\dim(W_1)$ and  $C_2=\max\{200 A \eta_1(1), 2 (\eta_1(1))^2\}$.  
   Then  there is a hyperplane $H$ in $\pi_1(L)$  passing through $\pi_1(q)$  and one component $H_-$ of $\pi_1(L)\backslash H$ such that $\text{l}_F(x)\le  (C_2)^{2m} \text{L}_F(q)$  for all $x\in L\cap  \pi_1^{-1}(H_-)$.

\end{Le}

\b{proof}
Let $S$ denote the unit tangent space  of $\pi_1(L)$ at $\pi_1(q)$.   We shall define two subsets  $G$, $B$ of $S$.  
 A point $s\in S$ lies in $G$ if  $\text{L}_F(x)\le 
 (C_2)^m \text{L}_F(q)$  for  every $x\in L$  
  such that $\pi_1(x)$ lies in the direction of $s$.
   A point $s\in S$ lies in $B$ if  $\text{l}_F(x)> 
 (C_2)^{2m} \text{L}_F(q)$  for    some  $x\in L$  such that
$\pi_1(x)$ lies in the direction of $s$.
  Clearly $G\cap B=\emptyset$.    Let $s_1\in S$ be the direction of $\pi_1(p)$,  and $s_2\in S$ the point in $S$ opposite to $s_1$.    Then $s_1\in B$ since  
$\text{l}_F(p)> (C_2)^{2m} \text{L}_F(q)$.  Lemma \ref{key} implies $\text{L}_F(x)\le  2(\eta_1(1))^2 \text{L}_F(q)$ for any point $x\in L$  such that $\pi_1(q)\in \pi_1(x)\pi_1(p)$.  Hence $s_2\in G$.

Let $H(B)\subset S$ be the convex hull of $B$ in the sphere $S$.  By   Caratheodory's theorem  
    for any $y\in H(B)$,
  there are $m$ points $x_1, \cdots, x_m\in B$ such that $y$ lies in the  spherical simplex  $\Delta_1$  spanned by 
   $x_1, \cdots, x_m$.    Let $\Delta_i$ be the spherical
   simplex spanned by $x_i, \cdots, x_m$.  
  Then there are $y_{i}\in \Delta_{i}$ ($1\le i\le m-1$)   with   $y_1=y$  such that  $y_i\in x_i y_{i+1}$.  
   Since $x_i\in B$, there exists a point $p_i\in L$ such that $\pi_1(p_i)$  lies  in the direction of $x_i$  with 
   $\text{l}_F(p_i)> 
 (C_2)^{2m} \text{L}_F(q)$. 
  Let $q_{m-1}\in L$ be a point 
  such that  $\pi_1(q_{m-1})$   lies 
in the direction of $y_{m-1}$   and  $\pi_1(q_{m-1})\in \pi_1(p_{m-1})\pi_1(p_m)$.
    Inductively, let $q_i\in L$ $(1\le i\le m-2$) be a   point 
 such that $\pi_1(q_i)$   lies in the direction of $y_i$ 
  and   $\pi_1(q_i)\in \pi_1(p_i)\pi_1(q_{i+1})$.  

 We claim $\text{L}_F(x)> C_2^{2m-1} \text{L}_F(q)$ for every $x\in L$ such that
   $\pi_1(x)\in   \pi_1(p_{m-1})\pi_1(p_m)$;   in particular,
   $\text{L}_F(q_{m-1})> C_2^{2m-1} \text{L}_F(q)$.
  Suppose   the claim does not  hold. 
      Then 
$\text{L}_F(x)\le  C_2^{2m-1} \text{L}_F(q)$ for   some 
 $x\in L$ satisfying 
   $\pi_1(x)\in   \pi_1(p_{m-1})\pi_1(p_m)$. 
  Since $\text{l}_F(p_{m-1})> (C_2)^{2m} \text{L}_F(q)$,  we have
 $  \text{l}_F(p_{m-1})>  C_2 \text{L}_F(x)\ge  C_1 \text{L}_F(x)$.  Now Lemma \ref{key} implies
   $$\text{L}_F(p_m)\le 2 (\eta_1(1))^2  \text{L}_F(x) \le   C_2^{2m} \text{L}_F(q), $$ contradicting  $\text{l}_F(p_m)> (C_2)^{2m} \text{L}_F(q)$.  By   considering $q_i\in L$
    satisfying   $\pi_1(q_i)\in  \pi_1(q_{i+1})\pi_1(p_i)$  and using Lemma \ref{key}  one inductively  proves  that 
 $\text{L}_F(q_{i})> C_2^{m+i} \text{L}_F(q)$.
   In particular,   $\text{L}_F(q_{1})> C_2^{m+1} \text{L}_F(q)$.
 Since $\pi_1(q_1)$   lies  in the direction of    $y_1=y$, 
    we see that   $y\notin G$.  
  Since $y\in H(B)$ is arbitrary, we have 
$H(B)\cap G=\emptyset$.

Now $H(B)$ is a  non-empty convex subset of the sphere $S$ and its complement is non-empty. 
  It follows that there is an  open hemisphere in its complement.  Hence there is an open hemisphere in  the complement of $B$.   Now the Lemma follows from the definition of $B$.

\end{proof}

\b{Le}\label{nonzero}
  Suppose $\dim(W)\ge 2$.  Then, 
  for any bounded subset  $X\subset \pi_1(L)$,  there exist 
  a  subset   $E\subset L\cap  \pi_1^{-1}(X)$  with   full measure   in 
   $L\cap  \pi_1^{-1}(X)$
 and 
two  constants    $M_1, M_2>0$ such that 
  $\text{L}_F(x)\ge M_1$ and $\text{l}_F(x)\le M_2$ for all $x\in E$.  

\end{Le}

\b{proof} 
We first show that there is some $M_1>0$ such that 
   $\text{L}_F(x)\ge M_1$  for all $x\in L\cap \pi_1^{-1}(X)$.  
 Suppose there is a   sequence of points  $x_i\in L\cap  \pi_1^{-1}(X)$ such that 
$\text{L}_F(x_i)\ra 0$.   Fix a point $p\in L$   with $\text{l}_F(p)>0$.
  Such a point always exists since $F|_L: L\ra F(L)$ is  quasiconformal and so
 by Pansu's theorem  is a.e. on $L$ Pansu differentiable and the Pansu differential is a.e. on $L$ non-singular. 
   For all sufficiently large $i$ we have
    $\text{l}_F(p)>C^{2m}_2 \text{L}_F(x_i)$,  where $C_2$ is the constant in Lemma \ref{halfspace}. 
      By Lemma \ref{halfspace},
  there is a hyperplane $H_i$  of $\pi_1(L)$   passing through $\pi_1(x_i)$  and a component $H_{i,-}$ of $\pi_1(L)\backslash H_i$   such that 
  $\text{l}_F(x)\le C_2^{2m} \text{L}_F(x_i)$ for all $x\in  L\cap  \pi_1^{-1}(H_{i,-})$. 
  Since  $X$ is bounded,  
the sequence $\pi_1(x_i)$ is bounded.  Hence    a subsequence 
  $H_{i_j,-}$  
of  the half spaces 
$H_{i,-}$  converges to an open half space $H_-$  of  $\pi_1(L)$.  
  Since every $v\in H_-$ lies in  $H_{i_j,-}$  for all sufficiently large $j$  and 
$\text{L}_F(x_i)\ra 0$, it follows that $\text{l}_F(x)=0$ for all 
 $x\in L$ such that $\pi_1(x)\in H_-$.  
    This means that the quasiconformal map $F|_L: L\ra F(L)$  has $zero$ differential on non-empty  open subsets, contradicting   Pansu's theorem.  

We  next show that $l_F$ is essentially bounded from above on $L\cap \pi_1^{-1}(X)$.
   After   pre-composing  and post-composing with left translations we may assume 
 $L=F(L)=W$.   Now we can write 
$L\cap \pi_1^{-1}(X)=X\oplus W_2\subset W_1\oplus W_2$.   
  Let  $B\subset W_2$ be a   bounded  non-empty open subset of $W_2$. 
  Let   $X_1\subset X$  be the subset of $X$ consisting of all 
 $x\in X$ such that   for   every  $b\in B$ either the Pansu differential $dF{(x, b)}$ does not exist   or     $dF{(x, b)}$   is not a graded isomorphism.  
  By Pansu's theorem and Fubini's theorem, 
 $X_1$ has zero measure in $X$.    Notice that 
$E:=L\cap \pi_1^{-1}(X\backslash X_1)$  has full measure in $L\cap \pi_1^{-1}(X)$.
We shall show that $l_F$ is bounded from above on 
   $E$.  
  First of all,  as a quasisymmetric map,    $F|_L: L\ra F(L)$  maps bounded sets to bounded sets. So 
 $F(X\oplus B)$ is bounded.   Now the first paragraph applied to $F^{-1}$ implies that   there is a positive lower bound   for $\text{L}_{F^{-1}}$  on $F(X\oplus B)$.  
   By (\ref{e3.001}) there is  some $M_1>0$ such that $l_F(q_1)\le M_1$ for any $q_1\in X\oplus B$.     Now let $q\in E$ be arbitrary.  Then by the definition of $E$ there is some $q_1\in X\oplus  B$ such that $\pi_1(q)=\pi_1(q_1)$ and the Pansu differential $dF({q_1})$ is a graded isomorphism. 
By Lemma \ref{add1} and   (\ref{e3.0})  
 we have
 $$l_F(q)\le C_1 L_F(q_1)\le C_1 \eta(1) l_F(q_1)\le C_1 \eta(1) M_1.$$

\end{proof}

\b{Le}\label{quasisimi}
For each left coset $L$  of $W$,  there is some constant $C_L>0$ such that 
   $F|_L $ is a $(C^{2m+3}_2, C_L)$-quasi-similarity,
 where $m=\dim(W)$ and $C_2$ is the constant in Lemma \ref{halfspace}.   

\end{Le}


\b{proof}
First   consider the case when $\text{dim}(W)=1$.   
By Pansu's  differentiability theorem and the theorem 
 in \cite{BKR} on absolute continuity on almost all curves,  we conclude that for a.e. left coset $L$ the following hold:
  $F$ is absolutely continuous on $L$,  at   a.e. point $q\in L$, the Pansu differential $dF(q)$ is a graded isomorphism.   Fix any such $L$. 
    We first observe that    to show 
 $F|_L $ is a quasi-similarity,   it suffices to show that 
$\text{l}_F(p)\le   C_1\cdot \text{L}_F(q)$  for any $p, q\in L$, where $C_1$ is the constant in Lemma \ref{key}.
Suppose there are two points $p, q\in L$ such that
  $\text{l}_F(p)> C_1\cdot \text{L}_F(q)$.     By Lemma \ref{key},  
   $\text{L}_F(x)\ra 0$ as $d(p, x)\ra \infty$ ($x\in L$).    By (\ref{e3.001}) we have 
  $\text{L}_{F^{-1}}(y)\ge \text{l}_{F^{-1}}(y) \ra \infty$ as $d(y, F(p))\ra \infty$ ($y\in F(L)$).
  However,    $\text{l}_F(p)> C_1 \text{L}_F(q)$  implies 
 $\text{l}_{F^{-1}}(F(q))> C_1 \cdot \text{L}_{F^{-1}}(F(p))$.    Applying  Lemma \ref{key} 
  to $F^{-1}$    we obtain
  $\text{L}_{F^{-1}}(y)\ra  0$, which is a contradiction.  
  Hence for a.e. left  coset $L$ of $W$,  
$F|_L $ is a $(C^{2}_2, C_L)$-quasi-similarity  for some $C_L$.  By continuity, this holds for every left coset  of $W$.

From now on we assume    $\text{dim}(W)\ge  2$.  In this case, both $F|_L$ and $F^{-1}|_{F(L)}$ have the following properties:
   (1)  absolutely continuous, 
 (2) Pansu  differentiable a.e.  and the Pansu  differential is a.e.  a graded isomorphism,
   (3)    absolutely continuous on almost all curves.    
  It follows that to show $F|_L$   is  a  $(C_2^{2m+3}, C_L)$ quasi-similarity
  for some $C_L>0$, it suffices to show that 
    there is    a  set of full measure $E\subset L$ such that 
    $\text{l}_F(x)\le  \eta(1)\eta_1(1) C_2^{2m+1} \text{L}_F(y)$  for all $x, y\in  E$.  
    We shall prove  this  by contradiction.  So suppose the above statement is not true. 
 Then   in particular  there are two points  $p, q\in L$ such that
   $ \text{l}_F(p)>\eta(1)  \eta_1(1)  C^{2m+1}_2 \text{L}_F(q)$.
  By Pansu's theorem, we may assume  $\text{l}_F(p)<\infty$.
   
 We first  observe that it suffices to show that there is a  constant $b_0\ge 1$ such that 
$\text{l}_F(x)\le b_0$  for all $x$ in a full measure subset of $L$:   
  the condition   $\text{l}_F(p)> \eta(1)  \eta_1(1)  C^{2m+1}_2 \text{L}_F(q)$
  implies that   $\text{l}_{F^{-1}}(F(q))> \eta(1)  \eta_1(1)  C^{2m+1}_2 \text{L}_{F^{-1}}(F(p))$.
  Then Lemma \ref{key} implies  that  there is some $y_0\in F(L)$ such that 
$$\text{L}_{F^{-1}}(y_0)<\frac{1}{{b_0 C_2^{2m}}\eta_1(1)}\min\{1,  \text{l}_{F^{-1}}(F(q))\} .$$  
     By Lemma \ref{halfspace},  
  there is a hyperplane   $H'$  in $\pi_1(F(L))$   passing trough $\pi_1(y_0)$  and a component $H'_-$  of $\pi_1(F(L))\backslash H'$ such that
     $ \text{l}_{F^{-1}}(y)<   \frac{1}{b_0\eta_1(1)}$  for all $y\in F(L)\cap  \pi_1^{-1}(H'_-)$. 
  Since $F^{-1}$ is Pansu differentiable a.e.,   $ \text{L}_{F^{-1}}(y)<   \frac{1}{b_0}$  for a.e.   $y\in F(L)\cap  \pi_1^{-1}(H'_-)$. 
 It follows that 
   $\text{l}_F(x)> b_0$ for a.e.   $x\in  L\cap  F^{-1}(\pi_1^{-1}(H'_-))$, contradicting  the assumption that  
$\text{l}_F(x)\le b_0$  for all $x$ in a full measure subset of $L$.

Since  we assume  there are two points  $p, q\in L$ such that
   $\text{l}_F(p)> \eta(1)  \eta_1(1)  C^{2m+1}_2 \text{L}_F(q)$,  
 by Lemma   \ref{halfspace},
    there is a hyperplane $H_1$ in $\pi_1(L)$  passing through $\pi_1(q)$  and one component $H_{1,-}$ of $\pi_1(L)\backslash H_1$ such that 
 $$\text{l}_F(x)\le  (C_2)^{2m} \text{L}_F(q)<\frac{1}{\eta(1)\eta_1(1) C_2}  \text{l}_F(p)$$  for all $x\in L\cap \pi_1^{-1}(H_{1,-})$. 
   Now Pansu's theorem and the quasisymmetry condition on $F$ imply 
$\text{L}_F(x)<\frac{1}{\eta_1(1)C_2}  \text{l}_F(p)$  for a.e.  
    $x\in L\cap  \pi_1^{-1}(H_{1,-})$. 
    Let $\tau:  \pi_1(L)\ra \pi_1(L)$ be the geodesic symmetry about $\pi_1(p)$,  that is,  for any $v\in \pi_1(L)$,  $\tau(v)$ is such that $\pi_1(p)$ is the 
   midpoint of $v\tau(v)$.    Now  Lemma \ref{key}  implies that   for a.e.   $y\in L\cap \pi_1^{-1}(\tau(H_{1,-}))$  we have 
    $$\text{L}_F(y)\le 2 (\eta_1(1))^2 \text{L}_F(\tau(y))\le  \frac{2\eta_1(1)}{C_2}l_F(p).$$   
     Hence  $L_F(x)$ is essentially bounded on  
$L\cap \left(\pi_1^{-1}(H_{1,-})\cup \pi_1^{-1}( \tau(H_{1,-}))\right)$. 
  If $\pi_1(p)\in H_1$, then we are done since   now $\text{L}_F(x)$ (hence $l_F(x)$) is bounded on the full measure subset 
 $L\cap \left(\pi_1^{-1}(H_{1,-})\cup \pi_1^{-1}( \tau(H_{1,-}))\right)$  of $L$.    So we assume $\pi_1(p)\notin H_1$.   
Let $B_1\subset \pi_1(L)$ be the part of a cylinder in $\pi_1(L)$ with center line  passing through  $\pi_1(p)$,   
    perpendicular to the hyperplane $H_1$  and bounded between $H_1$ and $\tau(H_1)$.     Then $B_1$ is bounded.  By Lemma \ref{nonzero}
      there are  numbers $M_1, M_2>0$
  and a set $E_1\subset  L\cap \pi_1^{-1}(B_1)$  with full measure in 
  $ L\cap \pi_1^{-1}(B_1)$  
  such that 
   $\text{L}_F(x)\ge M_1$ and $\text{l}_F(x)\le M_2$ for  every   $x\in E_1$.

  Since  we assume 
    $\text{l}_F(x)$ is not essentially bounded,   there is some hyperplane $\tilde{H}_1$ 
 of $\pi_1(L)$   parallel to $H_1$
 such that (1) $F$ is  Pansu  differentiable at some   $q_1\in E_1\cap \pi_1^{-1}(B_1\cap \tilde{H}_1)$; 
 (2)  there is some $p_1\in \pi_1^{-1}(\tilde{H}_1)$   with  $C^{2m+1}_2\cdot(\eta(1))^2 M_2< \text{l}_F(p_1)<\infty$.  Since $F$ is  Pansu  differentiable at 
 $q_1$,  we have $\text{L}_F(q_1)\le \eta(1) \text{l}_F(q_1)\le \eta(1) M_2$.    So
  $\text{l}_F(p_1)>C^{2m+1}_2\eta(1)\cdot   \text{L}_F(q_1)$.  Now,    as indicated above,    Lemma \ref{halfspace}
  and  Lemma \ref{key}    imply   that there
  is a hyperplane $H_2$ of $\pi_1(L)$  passing through $\pi_1(q_1)$ and a component $H_{2, -}$  of $\pi_1(L)\backslash H_2$ such that
   $\text{L}_F(x)$  is     essentially 
bounded    on 
$\pi_1^{-1}(H_{2, -})$    and $\pi_1^{-1}(\tau_1(H_{2, -}))$, where $\tau_1$ is the geodesic symmetry in   $\pi_1(L)$  about $\pi_1(p_1)$.  
  If $\pi_1(p_1)\in H_2$, then we are done as indicated above. So we assume $\pi_1(p_1)\notin H_2$. In this case,  $H_1$ and $H_2$ are not parallel.    We proceed inductively and eventually find $m$ hyperplanes $H_1$,  
 $H_2$,  $\cdots$,  $H_m$,  half spaces $H_{i, -}$  and points $p_0:=p, p_1, \cdots,  p_{m-1}$ with the following properties:\newline
  (1)  $\text{L}_F(x)$  (hence $l_F(x)$)   is   essentially bounded on  $\pi_1^{-1}(Q)$, where $Q$ is 
the union 
  $Q:=\cup_i H_{i, -}\cup \cup_i \tau_{i-1}(H_{i, -})$, where $\tau_{i-1}$ is the geodesic symmetry about the point $\pi_1(p_{i-1})$;\newline
  (2)   The complement of $Q$ in $\pi_1(L)$ is compact.\newline
  By Lemma \ref{nonzero},  $\text{l}_F(x)$ is  essentially  bounded on $L\backslash \pi_1^{-1}(Q)$. It follows that 
  $\text{l}_F(x)$  is  essentially   bounded on $L$,  and we are done.

\end{proof}

Now Theorem \ref{main1} follows from Lemma \ref{quasisimi}
  and   the proof of  Theorem  \ref{induct}.

\section{Examples}\label{examples}

In this Section we   apply our results to 3  examples.

We call a Carnot algebra quasisymmetrically rigid if the corresponding Carnot group is  quasisymmetrically rigid.
  All the examples here are quasisymmetrically rigid.  The first example satisfies $r_1(\mathcal N)=1$, and so is non-rigid by \cite{O}.  Similarly for  the Carnot algebras in the second example as soon as   $r_1(\mathcal N)=1$ for the factor $\mathcal N$.   The last example is rigid since $r_1(\tilde {\mathcal N})=5$.


\noindent
{\bf{First Example}}.    Let   $\mathcal{N}=V_1\oplus V_2\oplus V_3$ be a $3$-step   Carnot algebra defined as follows.
The first layer   $V_1$ has  basis $\{X, Y,  X', Y',  X_1, X_2\}$,  the second  layer $V_2$ has  basis $\{Z,  X_{1,2}\}$, and the third layer  $V_3$ has  basis 
  $\{X_{1,2,1}, X_{1,2,2}\}$.  The only non-trivial bracket relations are 
 $$[X, Y]=[X', Y']=Z,  \;  [X_1, X_2]=X_{1,2},\;  [X,  X_2]=[Y,  X_2]=X_{1,2},\; $$
  $$[X_{1,2}, X_1]=[X_{1,2},  X]=[X_{1,2},  Y]=X_{1,2,1},      [X_{1,2}, X_2]=X_{1,2,2}.$$
  Let $0\not=v\in V_1$.  Write $v=aX+bY+a_1 X_1+b_1X_2+a'X'+b'Y'$.     A direct calculation shows
$$[v,  X]=-b Z-b_1 X_{1,2},$$  
 $$[v,  Y]=aZ-b_1 X_{1,2},$$   
\b{equation}\label{e6.10}
[v,  X_1]=-b_1  X_{1,2},
\end{equation}   
$$[v,  X_2]=(a+b+a_1)X_{1,2},$$
 $$[v, X']=-b' Z,$$
  $$[v, Y']=a'Z,$$
$$[v, Z]=0,$$
   $$[v, X_{1,2}]=-b_1  X_{1,2,2}-(a+b+a_1) X_{1,2,1}. $$
  It   is easy to check  that 
  $\text{rank}(v)\ge 1$  and  
$\text{rank}(v)=1$  if and only if  $b_1=0$ and $a+b+a_1=0$.  
  Hence $W_1=\{aX+bY+a_1 X_1+a'X'+b'Y':  a+b+a_1=0\}$  is a non-trivial proper subspace  of  $V_1$  invariant under the action of   
    $\text{Aut}_g(\mathcal N)$.
    Notice that  $ X-Y, X-X_1, X', Y'$  form  a    basis for $W_1$   and 
  the subalgebra  
  $<W_1>=W_1\oplus \R Z$ 
       is isomorphic to the second Heisenberg algebra.  
   Since $N$ is $3$-step, one can not apply Theorem   \ref{main1}.   Also 
    $<W_1>$  is   a Heisenberg algebra   and so Theorem \ref{main3}  is not applicable
(there exist  non-biLipschitz quasiconformal maps on the Heisenberg groups  \cite{B}).  
  However, (\ref{e6.10})  implies 
 $[X_1, W_1]=0$.  
  Since   $X_1\notin W_1$,   by Theorem   \ref{main2},   
  $\mathcal N$ is quasisymmetrically rigid.


\noindent
{\bf{Second  Example.}}     Let $\mathcal N=V_1\oplus V_2$,  $\mathcal N'=V_1'\oplus V_2'$  be  $2$-step 
    Carnot algebras.   If   $r_1(\mathcal N)\not=r_1(\mathcal N')$,   then the discussion in Section
  \ref{product}   and Theorem \ref{main1}   show that the direct product $\mathcal N\oplus \mathcal N'$ and all their central products are 
 quasisymmetrically rigid.    If  $r_1(\mathcal N)\ge 1$  and  $r_1(\mathcal N')\ge 1$,  then Lemma \ref{level1} and Theorem \ref{main1} imply that all 
 level one products of $\mathcal N$ and $\mathcal N'$ are 
quasisymmetrically rigid.

\noindent
{\bf{Third Example}}.       Next we construct a Carnot algebra from the model Filiform algebra ${\mathcal F}^4$ and the free nilpoptent Lie algebra 
 ${\mathcal  F}_{2,6}$ ($2$-step with $6$ generators). 
  Recall that   ${\mathcal F}^4=V_1\oplus V_2\oplus V_3\oplus V_4$  with $V_1$ spanned by $e_1, e_2$,
  $V_2$ spanned by $e_3$, $V_3$ spanned by $e_4$ and $V_4$ spanned by $e_5$.  The only non-trivial bracket relations are 
 $[e_1, e_i]=e_{i+1}$ for $i=2,3,4$.  
  The   free nilpotent Lie algebra  ${\mathcal  F}_{2,6}=V'_1\oplus V'_2$, where $V'_1$ is spanned by $\{X_i:   1\le i\le 6\} $
  and $V'_2$ is spanned by $\{X_{ij}:  1\le i<j\le 6\}$.
 The only non-trivial bracket relations are $[X_i, X_j]=X_{ij}$  for $i<j$.   
  The Carnot algebra 
  $\mathcal  N$  is obtained from ${\mathcal F}^4\oplus  {\mathcal  F}_{2,6}$  by adding the following bracket relations:
  $$[e_2,   X_2]=X_{13},  \; [e_2,   X_3]=X_{24},  \; [e_2,   X_4]=X_{35}, $$
  $$ [e_2, X_6]=X_{15},\;  \; [e_1,   X_1]=X_{26}, \; [e_1, X_5]=X_{46}.$$
     Now it is easy to check that $\text{rank}(e_1)=\text{rank}(e_2)=5$.  Furthermore, a tedious calculation shows   
  $\text{rank}(v)\ge 6$    for  
$v\in (V_1\oplus V'_1)\backslash V_1$.  Hence 
   $W_1:=V_1$ is a non-trivial proper subspace of 
  $V_1\oplus V'_1$ invariant  under the action of  $\text{Aut}_g(\mathcal N)$
  and $<W_1> ={\mathcal F}^4$.   It was  proved in 
  \cite{X2}   that   ${\mathcal F}^4$  is  quasisymmetrically rigid.
  Now Theorem \ref{main3}    implies that   $\mathcal N$ is also quasisymmetrically rigid.


The group $N$ is $4$-step,   so   Theorem \ref{main1}   does not apply.    It also does not satisfy the condition in Theorem \ref{main2}. 
 Indeed, 
   let $X=a_1e_1+a_2e_2+\sum_{i=1}^6  b_iX_i$.  Now
  $[e_1,  X]=a_2e_3+b_1X_{26}+b_5X_{46}$ and $[e_2, X]=-a_1e_3+b_2X_{13}+b_3X_{24}+b_4X_{35}+b_6X_{15}$.
   So  $[X, W_1]\subset <W_1>$  
if and only if $X\in W_1$.

 \addcontentsline{toc}{subsection}{References}

\noindent Address:

\noindent Xiangdong Xie: Dept. of Mathematics  and   Statistics,   Bowling Green  State  University, 
  Bowling Green,  OH,   U.S.A.\hskip .4cm E-mail:   xiex@bgsu.edu

\end{document}